\documentclass[12pt,a4paper]{article}
\usepackage{amsmath,amssymb,amsthm}
\usepackage{showkeys}

\newtheorem{theorem}{Theorem}
\newtheorem{lemma}[theorem]{Lemma}
\newtheorem{example}[theorem]{Example}
\newtheorem{proposition}[theorem]{Proposition}

\newtheorem{corollary}[theorem]{Corollary}

\newtheorem*{theorema}{Theorem A}

\newtheorem*{hypothesis}{Hypothesis}

\newcommand{\stepone}{Step 1 (Bound)}
\newcommand{\steptwo}{Step 2 (Simplicity)}
\newcommand{\stepthree}{Step 3 (Irreducibles)}
\newcommand{\stepfour}{Step 4 (Reducibles)}
\newcommand{\stepfive}{Step 5 (Primitivity)}
\newcommand{\stepsix}{Step 6 (Imprimitivity)}

\newcommand{\Z}{\mathbb{Z}}

\newcommand{\Galph}{G_\alpha}
\newcommand{\Lalph}{L_\alpha}

\newcommand{\curlyc}[1]{\mathcal{C}_#1}
\newcommand{\curlycg}{\mathcal{C}(G)}

\newcommand{\point}{\Pi}
\newcommand{\lines}{\Lambda}
\newcommand{\linel}{\mathfrak{L}}
\newcommand{\hatt}{\hat {\ } }
\newcommand{\spaceS}{\mathcal{S}}

\newcommand{\sss}{ \subset  }

\def\itf{{It follows that }}

\newcommand{\up}{^{-1}}

\begin{document}
\pagestyle{plain} \pagenumbering{arabic}
\title{Large dimensional classical groups and linear spaces}
\author{Alan R. Camina, Nick Gill, A.E. Zalesski\thanks{The authors wish to thank the London Mathematical Society for financial support during the writing of this paper. The second author would also like to thank the University of Bristol for providing practical support.}}
\date{4 April 2007}\maketitle

\begin{abstract}
Suppose that a group $G$ has socle $L$ a simple large-rank
classical group. Suppose furthermore that $G$ acts transitively on
the set of lines of a linear space $\mathcal{S}$. We prove that,
provided $L$ has dimension at least $25$, then $G$ acts
transitively on the set of flags of $\mathcal{S}$ and hence the
action is known. For particular families of classical groups our
results hold for dimension smaller than $25$.

The group theoretic methods used to prove the result (described in
Section \ref{section:almostsimple}) are robust and general and are
likely to have wider application in the study of almost simple
groups acting on finite linear spaces.
\end{abstract}

\section{Introduction}
A {\it linear space} $\spaceS$ is an incidence structure
consisting of a set of points $\point$ and a set of lines $\lines$
in the power set of $\point$ such that any two points are incident
with exactly one line. The linear space is called {\it
non-trivial} if every line contains at least three points and
there are at least two lines. The linear space is called {\it
finite} provided $\point$ and $\lines$ are finite. All linear
spaces which we consider are finite; we write $v=|\point|$ and
$b=|\lines|$.

 This paper is part of a sequence attempting to classify those groups which can act line-transitively on a finite linear space. In \cite{campraeger, campraeg2} it was shown that such a group takes one of three forms: Either it contains a normal subgroup which acts intransitively on the set of points, or it contains an elementary-abelian normal subgroup acting regularly on the set of points, or it is almost simple. This last case is equivalent to the group having {\it simple socle} (recall that the socle of a finite group $G$ is the product of the minimal normal subgroups of $G$).

 The results in \cite{campraeger, campraeg2} have inspired the study of almost simple groups acting line-transitively. In particular studies have been made when the socle is a sporadic group (\cite{camspiez}), an alternating group (\cite{cnp}) or a member of particular families of low rank groups of Lie type (\cite{gill1, liupsl2,
liusuzuki,liug2even, liug2odd, liupsu3even, liu3d4, liuree,
liuree2,  Zhouliu,Zhouree, zhouree2}). We continue this
investigation by considering the situation when the socle of a
line-transitive automorphism group $G$ is $L$ a simple large-rank
classical group.

Our major result is that, in this situation, any line-transitive
action is in fact flag-transitive. (Here a {\it flag} is an
incident point-line pair.) The flag-transitive actions of almost
simple groups were fully classified in \cite{ftclass} (with proofs
in \cite{bdd2, davies, delandt4, delandt3, kleidman4, saxl, fed}) and
are well-known. When the socle is a large-rank classical group
there is one infinite family of flag-transitive actions:

\begin{example}\label{example:flagtransitive}
We have $\spaceS=PG(n-1,q)$, projective space of dimension
$n-1\geq 2$ over a field of size $q$ where $q$ is a prime power.
Any group $G$ with $PSL_n(q)\leq G\leq P\Gamma L_n(q)$ acts
flag-transitively on $\spaceS$.
\end{example}

In order to state our theorem we need to define a triple $(N_1,
N_2,N_3)\in\Z^3$, $N_1\leq N_2\leq N_3$, dependent on the socle
$L$. Values for this triple are given in Table \ref{table:triple}.

\begin{table}
\begin{center}
\begin{tabular}{|c|c|c|c|c|}
\hline \hline
$L$ & Conditions & $N_1$ & $N_2$ & $N_3$\\
\hline \hline
$\Omega_n(q)$ & $n$ odd, $q$ odd & 7 & 13 &21\\
\hline
$P\Omega^\epsilon_n(q)$ & $n$ even, $q$ odd & 18 & 26& 26\\
\hline
$\Omega^\epsilon_n(q)$ & $n$ even, $q$ even & 16 & 26& 26\\
\hline
$PSU_n(q)$ & $nq$ odd & 11 & 15& 15\\
\hline
$PSU_n(q)$ & $n$ even, $q$ odd & 16 & 22& 22 \\
\hline
$PSU_n(q)$ & $q$ even & 11 & 13& 13\\
\hline
$PSp_n(q)$ & $n$ even, $q$ odd & 12 & 22& 22\\
\hline
$Sp_n(q)$ & $n$ even, $q$ even & 8 & 14& 14\\
\hline
$PSL_n(q)$ & $nq$ odd & 7 & 17 & 17\\
\hline
$PSL_n(q)$ & $n$ even, $q$ odd & 12 & 22& 22\\
\hline
$PSL_n(q)$ & $q$ even & 8 & 17 & 17\\
\hline
\end{tabular}
\caption{Values for $(N_1,N_2,N_3)$}\label{table:triple}
\end{center}
\end{table}

\begin{theorema}
Let $G$ be a group which acts transitively on the set of lines of
a linear space $\spaceS$. Suppose that $G$ has socle $L$ a simple
classical group of dimension $n$. The following statements hold:
\begin{itemize}
\item   If $n\geq N_1$ then $L$ acts transitively on the set of
lines of $\spaceS$.
\item If $n\geq N_2$ and $L$ acts primitively on the set of points of $\spaceS$ then we have
Example \ref{example:flagtransitive}.
\item If $n\geq N_3$ then we have
Example \ref{example:flagtransitive}.
\end{itemize}
\end{theorema}

\begin{corollary}\label{corollary:first}
Let $G$ be a group which acts transitively on the set of lines of
a linear space $\spaceS$. Suppose that $G$ has socle $L$ a simple
classical group of dimension $n\geq 25$. Then $G$ acts transitively on
the set of flags of $\spaceS$ and we have Example
\ref{example:flagtransitive}.
\end{corollary}

The rest of the paper will be occupied with proving Theorem A. In Section \ref{section:background} we will outline a number of general background lemmas concerning groups acting line-transitively on finite linear spaces; in particular we record several lemmas which first appeared in \cite{cnp} and which will be crucial to our proof of Theorem A. In Section \ref{section:almostsimple} we outline a method for applying these results to the case where $G$ is almost-simple. Then, in the remaining sections, we apply the method of Section \ref{section:almostsimple} to the different families of almost simple groups listed in Theorem A. Note that in our final section we are able to state a stronger result than Corollary \ref{corollary:first} by referring to some cases not covered by Theorem A.

We write $\alpha$ for a point of $\spaceS$, $\linel$ for a line of $\spaceS$ and $G_\alpha$, $G_\linel$ the respective stabilizers in $G$.

\section{Some background lemmas}\label{section:background}

We list here some well-known lemmas which we will use later.
Suppose throughout this section that $G$ acts line-transitively on
a linear space $\spaceS$. Block\cite{block} proved that
line-transitivity implies point-transitivity. Thus, if a linear
space $\spaceS$ is line-transitive then every line has the same
number, $k$, of points and every point lies on the same number,
$r$, of lines. Such a linear space is called {\it regular}; all
linear spaces which we consider from here on will be assumed to be
regular. The first lemma is proved easily by counting and holds
for all regular finite linear spaces.

\begin{lemma}\label{lemma:arith}
\begin{enumerate}
\item   $b=\frac{v(v-1)}{k(k-1)}\geq v$ (Fisher's inequality);
\item   $r=\frac{v-1}{k-1}\geq k$;
\end{enumerate}
\end{lemma}

\begin{lemma}\label{lemma:davies}\cite{davies, camsiem}
If $g\in G$ then $g$ fixes at most $r+k-3$ points.
\end{lemma}

\begin{lemma}\cite[Lemma2]{camsiem}\label{lemma:invfixespoints}
If $g\in G$ is an involution and $g$ fixes no points then $G$ acts
flag-transitively on $\spaceS$.
\end{lemma}

If $\spaceS$ is not a projective plane then by Fisher's inequality
$b>v$ and, since $b=v(v-1)/(k(k-1))$, there must be some prime $t$
that divides both $v-1$ and $b$. We shall refer to such a prime as
a {\it significant} prime of $\spaceS$.

\begin{lemma}\label{lemma:sigprime}\cite[Lemma 6.1]{cnp}
Suppose that $\spaceS$ is not a projective plane and let $t$ be a
significant prime of $\spaceS$. Let $S$ be a Sylow $t$-subgroup of
$G_\alpha$. Then $S$ is a Sylow $t$-subgroup of $G$ and $G_\alpha$
contains the normalizer $N_G(S)$.
\end{lemma}

\begin{lemma}\label{lemma:selfnormalizing}
$N_G(G_\alpha)=G_\alpha$.
\end{lemma}
\begin{proof}
Let $S\in Syl_t G_\alpha$ and suppose that $N_G(G_\alpha)>
G_\alpha$. Then $N_G(S)>N_{G_\alpha}(S)$ hence $\Galph$ does not
contain the normalizer of a Sylow $t$-subgroup for any prime $t$.
This contradicts Lemma \ref{lemma:sigprime}.
\end{proof}


Finally we present a series of results which first appeared in
\cite{cnp} and which will be central to our analysis of the almost
simple case.



\begin{lemma}\cite[Lemma 2.2]{cnp}\label{one}
Let $g\in G$ fix $f$ points. Then $k<2v/f$.\end{lemma}




\begin{lemma}\cite[Corollary 2.4]{cnp}\label{bound1}
Let $g \in G_\alpha$. Assume that $g$ has $w$ conjugates in $G$
and that $a$ of them lie in $G_\alpha$. Then $$k<\frac{2w}{a}.$$
\end{lemma}



\begin{lemma}\cite[Proposition 2.7]{cnp}\label{lemma:main}
Let $h,g\in G$ with $g\in\Galph$ and $h\not\in\Galph$. Define
$$w=|G:C_G(g)|, \ c = |G:C_G(h)|\geq |\Galph:C_{\Galph}(h)|. $$
and suppose that $g$ has $a$ conjugates in $\Galph$. Then $r\leq
kc$, $k<\frac{2w}{a}$ and $$v=r(k-1)+1\leq
ck(k-1)+1<\frac{2wc}{a}(\frac{2w}{a}-1)+1\leq\frac{4w^2c}{a^2}.$$
\end{lemma}
\begin{proof}
Look at $K=G_\alpha\cap G_{h\alpha}$. This subgroup has to lie in
the stabilizer $G_{\linel}$
 of some line $\linel$. Then $$|G:K|=|G:G_{\linel}|\cdot|G_{\linel}:K|=|G:G_\alpha|\cdot|G_\alpha:K|.$$

\itf
$$|G_\alpha:K|=\frac{|G:G_{\linel}|\cdot|G_{\linel}:K|}{|G:G_\alpha|}=\frac{b|G_{\linel}:K|}{v}
~{\rm so}~ r\leq k|G_\alpha:K|.$$

We need to obtain an upper bound for the right hand side. Let
$C=C_G(h)$ and $L=G_\alpha\cap C$. Then $[h,L]=1$ so $hLh\up
=L\sss G_{h\alpha}$ so $L\sss G_\alpha\cap G_{h\alpha}=K.$ Hence
$|G_\alpha:K|\leq |G_\alpha:L|$. Set $c=|G:C|$ and consider
$\Omega =G/C$ as a $G$-set. Let $x\in  \Omega $ be the point
corresponding to the coset $C$. Then the stabilizer in $G_\alpha$
of $x$ is $G_\alpha\cap C=L$. Therefore, the orbit $G_\alpha x$ is
of size $|G_\alpha:L|$. Therefore, $|G_\alpha:L|\leq | \Omega |=c$
so $|G_\alpha:K|\leq c$. This implies the first inequality in the
lemma. The second follows from Lemma \ref{bound1}:
$$v=r(k-1)+1\leq |G_\alpha:K|\cdot ck(k-1)+1<
\frac{2w}{a}\left(\frac{2w}{a}-1\right)+1.$$
\end{proof}

We can slightly strengthen this result in the following ways:

\begin{corollary}
\begin{itemize}
\item If $h$ is an involution then $r\leq \frac{kc}{2}$ and
$$v<\frac{wc}{a}(\frac{2w}{a}-1)+1<\frac{2w^2c}{a^2}.$$

\item In general $|G:C_{\Galph}(g)|\leq 4w^2c.$
\end{itemize}
\end{corollary}
\begin{proof}
To prove the first statement suppose that $h$ is an involution.
Then, in the notation of Lemma \ref{lemma:main}, $h\in
G_\linel\backslash K$. This implies that

$$|G_\alpha:K|=\frac{|G:G_{\linel}|\cdot|G_{\linel}:K|}{|G:G_\alpha|}\geq\frac{2b}{v}
~{\rm so}~ r\leq \frac12 k|G_\alpha:K|.$$

To prove the second statement observe that, for $g\in\Galph$, it
is clear that $a\geq|\Galph:C_{\Galph}(g)|$. This implies that
\begin{eqnarray*}
&&|G:\Galph|=v<\frac{4w^2c}{a^2}\leq \frac{4w^2c}{|\Galph:C_{\Galph}(g)|^2} \\
&\implies& |G:C_{\Galph}(g)|\leq
\frac{4w^2c}{|\Galph:C_{\Galph}(g)|}\leq 4w^2c.
\end{eqnarray*}
\end{proof}








\section{The almost simple case}\label{section:almostsimple}

Lemma \ref{lemma:main} will be our primary tool in the analysis of
an almost simple group $G$ acting transitively on the set of lines
of a linear space $\spaceS$. We first outline our notation for the
rest of the paper.

\subsection{Notation}

Let $L$ be the socle of $G$. We consider $L$ in the different
families of classical groups; thus $L$ is a simple group of
dimension $n$ defined over a field of size $q$ where $q=p^a$ for
some prime $p$. Let $\curlyc{i}, i=1,\dots, 8$ be the Aschbacher
families of maximal subgroups of $G$ as described in \cite{kl}. The
collection of all maximal subgroups in these families will be
written $\curlycg$. Let $L^*$ be the linear group which covers
$L$. For $H<L$ (resp. $g\in L$) let $H^*$ (resp. $g^*$) be the
pre-image of $H$ (resp. $g$) in $L^*$ under the natural
homomorphism.

Take $V$ to be the $n$-dimensional vector space over the field of
size $q$ on which $L^*$ acts naturally. Let $\kappa$ be the
non-degenerate form on $V$ which is preserved by $L^*$, i.e.
$L^*\unlhd I(V,\kappa)$, the set of isometries of $\kappa$ in
$GL(V)$. We take $S(V,\kappa)$ to be the set of isometries of
$\kappa$ in $SL(V)$ while $\Omega(V,\kappa)$ coincides with
$S(V,\kappa)$ except in the orthogonal case when it is of index
$2$ in $S(V,\kappa)$. We will often write simply $I(V), S(V)$ or
$\Omega(V)$ when the form of $\kappa$ is clear.

We will reserve lower case Greek letters to represent the symbols
$+$ and $-$; in particular $\Omega_n^\zeta(q)$ is used to
represent one of the groups $\Omega_n^+(q)$ or $\Omega_n^-(q)$.

For a group $G$ take $P(G)$ to be the degree of the minimum
permutation representation of $G$. We represent cyclic groups of
order $c$ by the integer $c$, while soluble groups of order $c$
will be represented by $[c]$. An extension of a group $G$ by a
group $H$ will be written $G.H$; if the extension is split we
write $G:H$. We will sometimes precede the structure of a subgroup
$H$ of a projective group with $\hatt$ which means that we are
giving the structure of the pre-image in the corresponding
universal group (i.e. $H^*$). The notation $\frac12 G$ refers to a
normal subgroup in $G$ of index $2$; $\mathit{O}_p(G)$ refers to
the largest normal $p$-group in $G$. We write $|H|_p$ for the
highest divisor of $|H|$ which is a power of a prime $p$;
similarly $|H|_{p'}=\frac{|H|}{|H|_p}$.

\subsection{Our Method}

The remainder of the paper will be concerned with proving Theorem
A. We operate from now on under the following hypothesis.

\begin{hypothesis}
Let $L\unlhd G\leq Aut L$ where $L$ is a finite simple classical
group defined over a finite field with $q=p^a$ elements, for a
prime $p$. Suppose that $G$ acts line-transitively but not
flag-transitively on a linear space $\spaceS$ with parameters
$(b,v,r,k)$. Suppose that $\Lalph\leq M<L$.

Let $n$ be the dimension of the classical geometry for $L$. For
each family of finite simple groups which we consider we will
define three integers $N_1\leq N_2\leq N_3$.
\begin{itemize}
\item   For Steps 1 and 2 we assume that $n\geq N_1$.
\item From Step 3 (Irreducibles) on we assume that $n\geq N_2$ and that $G=L$.
\item For Step 6 (Imprimitivity)  we assume that $n\geq N_3$.
\end{itemize}

\end{hypothesis}

We will prove in Steps 1 and 2 that, for $n\geq N_1$, $L$ acts
line-transitively on $\spaceS$. Then we will prove that, for
$n\geq N_2$ and $G$ acting point-primitively, no actions exist;
then in Step 6 (Imprimitivity) we prove that, for $n\geq N_3$, no
actions exist.

Two particular situations can be excluded immediately.

First of all \cite[Theorem A]{gill2} implies that $\spaceS$ is not
a projective plane. In particular we must have $b>v$ and so
$\spaceS$ has a significant prime.

Secondly \cite[Theorem 1]{gill3} implies that the characteristic
prime, $p$, of $L$ is not significant.

We now describe the steps that we will use to prove Theorem A in
each case. In the course of our description we prove a number of
lemmas which can be seen to hold for a much weaker hypothesis; in
particular the lower bounds on dimension are often not used.

\begin{itemize}
\item {\bf \stepone}: Find small values for $w$ and $c$ (we use
\cite{kl, gorenstein}). Note that all involutions must fix a point
so we can take $g$ to be an involution (Lemma
\ref{lemma:invfixespoints}). Use Lemma \ref{lemma:main} to get a
rough bound for $k$ and $v$.

\item {\bf \steptwo}: We apply the principle of ``exceptionality":
Let $B$ be a normal subgroup in a group $G$ which acts upon a set
$\Pi$. Then $(G,B,\Pi)$ is called {\it exceptional} if the only
common orbital of $B$ and $G$ in their action upon $\Pi$ is the
diagonal (see \cite{gms}). We refer to the following lemma:

\begin{lemma}\cite[Lemma 26]{gill1}\label{lemma:exceptional}
Suppose a group $G$ acts line-transitively on a linear space
$\spaceS$. Let $B$ be normal in $G$ such that $|G:B|$ is a prime.
If $B$ is not line-transitive on $\spaceS$ then either $\spaceS$
is a projective plane or $(G,B,\Pi)$ is exceptional.
\end{lemma}

We know that $\spaceS$ is not a projective plane. Suppose that $L$
is not line-transitive on $\spaceS$. Then there exist groups
$G_1,G_2$ such that $L\unlhd G_1\lhd G_2\leq G\leq Aut L$ where
$|G_2:G_1|$ is a prime and $G_1$ is not line-transitive on
$\spaceS$ while $G_2$ is. By  Lemma \ref{lemma:exceptional},
$(G_2, G_1, \Pi)$ is an exceptional triple and \cite[Theorem
  1.5]{gms} implies that $\Lalph$ lies inside $M$ a maximal subgroup of $G$ and $M\in\curlyc{5}$.
In fact in all cases we know that $|M|\leq |\hatt G(q^\frac15)|$
where $G(q^\frac15)$ is the group of similarities of the same
classical geometry over a field of size $q^\frac15$. We will use
this to yield a contradiction for $n\geq N_1$.

For the remaining steps we will assume that $G=L$, that is that
$G$ is simple. This will often enable us to refine our estimate of $w$ and $c$ to improve the strength of our upper bound on $v$.

\item {\bf \stepthree}: We prove that, for $n\geq N_2$, $\Lalph$ must lie in a reducible subgroup. Our analysis
primarily makes use of the theorems provided by Liebeck in
\cite{liebeck3}. Liebeck lists the irreducible subgroups $M$ with
minimal index; it is enough for us to show that $|L:M|$ is greater
than the upper bound for $v$.

When $L=PSp_n(q)$ with $q$ even the results of \cite{liebeck3} are
not quite strong enough and we will also need to make use of the
following lemma:

\begin{lemma}\label{lemma:notcurly}\cite[Theorem 5.2.4]{kl}
Let $L$ be a classical simple group with associated geometry of
dimension $n$ over a field of size $q$. Also let $L\unlhd G\leq
Aut(L)$, and let $H$ be any subgroup of $G$ not containing $L$.
Then either $H$ is contained in a member of $\curlycg$ or one of
the following holds:
\begin{enumerate}
\item $H$ is $A_m$ or $S_m$ with $n+1\leq m\leq n+2$; \item
$|H|<q^{3n}$.
\end{enumerate}
\end{lemma}

\item{\bf \stepfour}: We (usually) have two possibilities for
maximal reducible subgroups - those which are parabolic and those
which stabilize a non-degenerate subspace. We are able to rule out
the first case in almost all situations as follows.

\begin{lemma}\label{lemma:parabolics}
Suppose that $L$ is a finite simple group isomorphic to one of
$PSp_n(q)$ ($n\geq 4$), $PSU_n(q)$ ($n>2$), $\Omega_n(q)$ ($n$
odd, $n\geq 7$), $P\Omega_{2s}^\epsilon(q)$ ($s\geq 4$). If
$\Lalph$ lies in a parabolic subgroup, $P$, of $L$ then $L =
P\Omega_{2s}^+(q)$ with $s$ odd. Furthermore $P=P_s$, the
stabilizer of an $s$-dimensional totally singular subspace.
\end{lemma}
\begin{proof}
Suppose that $\Lalph$ does lie in a parabolic subgroup $P$ of $L$.
Let $P_L$ be a Levi-complement of $P$. Now $P^*$ fixes $W$, a totally singular (or totally isotropic) subspace of
$V$. We can take a basis for $W$, $\{e_1,\dots, e_m\},$ such that there
exists $f_1, \dots f_m$ in $V$ with $(e_i, f_i)$ orthogonal hyperbolic
pairs. Furthermore we can choose $f_1,\dots f_m$ such that $P_L$ stabilizes
$\langle f_1,\dots f_m\rangle$.

We first show that there exists an element $g\in L^*$ such that $\langle e_i\rangle g=\langle
f_i\rangle$ and $\langle f_i\rangle g = \langle e_i\rangle$ for $i=1,\dots,
m$.
\begin{itemize}
\item   If $(V,\kappa)$ is symplectic or unitary, then we can define $g$ to
be such that $e_ig=f_i$ and $f_ig=-e_i$ for $i=1,\ldots ,m$, with
the trivial action on the orthogonal complement of $\langle e_1,
\ldots e_m, f_1\ldots f_m\rangle$.

\item   Let $(V,\kappa)$ be orthogonal and take $X$ to be the orthogonal
complement of all $e_i,f_i$. Let $\sigma$ be the projection of
$I(V,\kappa)$ into $I(V,\kappa)/L^*$ which is a group of exponent
2 and of order $2(2,q-1)$. Define $g_i\in I(V,\kappa)$ as follows:
$e_ig_i=f_i$,  $f_ig_i=e_i$ and $vg_i=v$ for all $v\in(e_i^\perp
\cap f_i^\perp)$. Set $g=g_1\cdots g_m$.

By Witt's theorem, all $g_i$ are conjugate in $I(V,\kappa)$, hence
$\sigma(g_i)=\sigma (g_j)$ for $1\leq i,j\leq m$. If $m$ is even
then $\sigma (g)=\sigma(g_1)\cdots \sigma(g_m)=1$ as
$I(V,\kappa)/L^*$ is of exponent 2. So $g\in L^*$.

Suppose $m$ is odd. If $\dim X>1$ then there exists an element
$y\in I(V,\kappa)$ fixing all $e_i,f_i$ such that $\sigma
(y)=\sigma (g_m)$. Then $y\in P^*$ and $gy\in L^*$. So the claim
follows by replacing $g$ with $gy$.

The case $\dim X=0$ corresponds to the possibility given in the
theorem, as $X=0$ means that $I(V,\kappa)=O^+(V)$.

Finally let $\dim X=1$; then $n=2m+1$ and $q$ is odd. Set
$T=\langle e_m,f_m,X\rangle$. Let $\mu : I(T,
\kappa|_T)\rightarrow I(V,\kappa)$ be
 the natural embedding.  It is well known that $\mu (\Omega (T))\subset \Omega (V)$.

Now $PSL(2,q)\cong \Omega (T)$ and this  isomorphism is obtained from the adjoint action of $SL(2,q)$ on
the space $S$ of $(2\times 2)$-matrices over $F$ with zero trace on which
the symmetric bilinear form is defined as trace$(ab)$ for $a,b\in S$.
Consider the matrices

$$e=\begin{pmatrix}
0&1\\
0&0
\end{pmatrix},
~~~ f=\begin{pmatrix}
0&0\\
1&0
\end{pmatrix}
~~~t=\begin{pmatrix}
1&0\\
0&-1
\end{pmatrix}
~~~ h=\begin{pmatrix}
0&1\\
-1&0
\end{pmatrix}.$$

 Then
 $e,f,t\in S$ and $heh^{-1}=f$,  $hfh^{-1}=e$ and $hth^{-1}=-t$. Let $g_m'$
 be the image of $h$ in $ \Omega (T)$. Then we replace $g$ by $g'=g_1\cdots
 g_{m-1}\cdot \mu (g_m')$. The images of $e,f$ in $T$ can be taken for
$e_m,f_m$. It is easy to observe that then $g'\in P^*$ acts as $g$ on
$e_i,f_i$ and $g'\in \Omega (V)$.
\end{itemize}

Thus we have our element $g$ in all cases. Let $h$ be the
projection of $g$ in $L$. Then $h$ normalizes $P_L$ but $h$ does
not lie in $P_L$. Let $S\in Syl_t P_L$. Then $N_L(S)>N_P(S)$ hence
the only $t$ for which $P$ contains the normalizer of a Sylow
$t$-subgroup is for $t$ the characteristic prime. But this means
that $p$ is a significant prime for $\spaceS$ and this possibility
is excluded by \cite[Theorem 1]{gill3}.
\end{proof}

For Step 4 (Reducibles) we deal with those subgroups not excluded
by Lemma \ref{lemma:parabolics}. We apply the rough bound to those
maximal subgroups $M$ in $\curlyc{1}$ to give bounds on the
dimension of subspaces stabilized by $M$.

\item{\bf \stepfive}:
Although we have labelled this step ``Primitivity", the situations
we consider are in fact more general than the primitive case. Nonetheless, at the completion of Step 5 (Primitivity), Theorem A will be proved for the situation where $G$ acts primitively on the set of points of $\spaceS$. There are two situations which we need to consider separately:
\begin{itemize}
\item $L\neq PSL_n(q)$: In this case the primary remaining case is
when $\Lalph$ preserves a non-degenerate subspace of $V$.  We will
typically take $M$ to be the projective image of $M^*$, the
stabilizer of a non-degenerate $m$-dimensional subspace $W$ of
$V$. Then $M^*=(X_m(q) \times Y_{n-m}(q)).[s]$  where $[s]$ is a
small soluble group, $$(X,Y)\in\{(Sp,Sp), (SU,SU), (\Omega ^\zeta,
\Omega^\eta), (\Omega, \Omega^\eta), (\Omega ^\eta, \Omega)\}$$
and $m\leq\frac{n}{2}$. We set $U=W^\perp$ and write $I(W)$ for
the group of isometries of $(W,\kappa|_W)$, similarly for $I(U)$.
Formally, $$I(W)=\{g\oplus 1_U: g\in I(W,\kappa|_W)\}.$$

In this step we add one supposition to our hypothesis as follows: {\bf Suppose that $\Lalph^*\geq
\Omega(U)$, the quasi-simple subgroup normal in $I(U)$.} We prove that, with this extra supposition, our hypothesis leads to a contradiction. Note that $\Omega(U)\cong Y_{n-m}(q)$. We will make use of the following lemma:

\begin{lemma}\label{lemma:noncentral}
Suppose $q$ is odd and $(V,\kappa)$ is symplectic, orthogonal or
unitary. Take $\Omega(U)\leq\Lalph^*\leq M^*$. Then there exists
an involution $g$ in $\Lalph^*$ which is not central in
$\Omega(U)$. Furthermore we can choose $g$ so that it acts as the
identity on a $1$ (resp. $2$) dimensional subspace, $X$, of $V$
for $n$ odd (resp. $n$ even) while taking all vectors in $X^\perp$
to their negative.
\end{lemma}
\begin{proof}
Take $n$ odd. There exists a non-singular vector, $u$, in $U$.
Define $g$ to be the involution in $I(V)$ which fixes $u$ and
which takes all vectors in $\langle u \rangle ^\perp$ to their
negative. In the symplectic and unitary cases all such $g$ must
lie in $\Omega(V)$. In the orthogonal case this is not necessarily
true however if we choose $u$ carefully we can ensure that $g$
lies in $\Omega(V)$.

Clearly $g$ normalizes $\Omega(U)$ and centralizes $I(W)$. Thus
$g$ normalizes $\Lalph^*$ and hence, by Lemma
\ref{lemma:selfnormalizing}, lies in $\Lalph^*$. It is clear that
$g$ is not centralized by $\Omega(U)$.

When $n$ is even we do the same but instead of using a
non-singular vector in $U$ we use a hyperbolic pair or, in the
orthogonal case, we may take $U$ to be non-degenerate and
anisotropic.
\end{proof}

This covers most of the {\it point-primitive} situations. There
are some occasional other possibilities (in even characteristic,
or thrown up by Lemma \ref{lemma:parabolics}) which we also
consider in this step.

\item $L=PSL_n(q)$: In this case Lemma \ref{lemma:parabolics} does
not apply. Hence we need to consider the situation where
$\Lalph\leq P_m$ a parabolic subgroup associated with an
$m$-dimensional subspace of $V$.

In this step we add one supposition to our hypothesis as follows: {\bf Suppose that $H\leq\Lalph\leq P_m$ where $H\cong \hatt
SL_{n-m}(q)$ and $H$ is normal in a Levi-complement of $P_m$.} We prove that, with the addition of this supposition, our hypothesis implies that we have Example \ref{example:flagtransitive}. This also excludes all point-primitive situations.

Note that a similar line of argument to that given in Lemma
\ref{lemma:parabolics} allows us to conclude that $m< \frac{n}2$
(i.e. $m\neq \frac{n}2$).

\end{itemize}

\item{\bf \stepsix}: We are left with the {\it point-imprimitive}
situation. In this case we have the following result of Delandtsheer and Doyen:
\begin{proposition}\cite{deldoy}\label{prop:deldoy}
Suppose that $L$ acts transitively on the set of lines of a linear space $\spaceS$ and acts imprimitively on the set of points of $\spaceS$. Then, if $\Lalph$ is the stabilizer of a point, we have that
$$|L:M|<{k\choose2}, \ |M:\Lalph|<{k\choose2}$$
where $M$ is any group such that $\Lalph<M<L$.
\end{proposition}

When using this result we will most often take $M$ to be a maximal subgroup of $L$. Again we distinguish between when $L=PSL_n(q)$ and when $L\neq
PSL_n(q)$. The former case is reasonably straightforward so, for now, we consider only the latter case.

We suppose that $\Lalph$ stabilizes a non-degenerate subspace $W$
of dimension $m\leq\frac{n}2$. As before write $U=W^\perp$,
$\Omega(W)\cong X_m(q)$ etc. The following lemmas will be useful:

\begin{lemma}\label{lemma:imprim_parabolics}
Suppose $\Lalph<M=\hatt (X_m(q) \times Y_{n-m}(q)).[s]$ and
$\Lalph^*\cap\Omega(U)$ lies inside a maximal parabolic subgroup,
$P^*,$ of $\Omega(U)\cong Y_{n-m}(q)$. Then $m=1,
Y_{n-m}=\Omega_{n-1}^+$, $P=P_{\frac{n-1}2}$ and $q(\frac{n-1}2)$ is
odd.
\end{lemma}
\begin{proof}
Our proof is very similar to the proof of Lemma
\ref{lemma:parabolics}. \cite[Lemma 4.1.1]{kl} implies that
$$\Lalph^*\leq\Omega(W)(P^*.[x])$$ where $[x]$ is a subgroup of
$[s]$ and $P^*.[x]$ is isomorphic to a parabolic subgroup of some
automorphism group of $\Omega(U)$.

We know that, except in the case listed, there exists $$h\in
\Omega(U)\backslash P^*$$ which normalizes the Levi-complement,
$P_{L}^*$, of $P^*$. In fact it is not hard to see that it must
normalize $\Omega(W)(P_{L}^*.[x])$.

Now, as in Lemma \ref{lemma:parabolics}, let $S\in Syl_t
(P_L.[x])$. Then $N_L(S)>N_P(S)$ hence the only $t$ for which $P$
contains the normalizer of a Sylow $t$-subgroup is for $t$ the
characteristic prime. But this means that $p$ is a significant
prime for $\spaceS$ and this possibility is excluded by
\cite[Theorem 1]{gill3}.
\end{proof}

\begin{lemma}\label{lemma:imprim_nondeg}
Take $H<L^*$ and suppose that $H<M^*=(\Omega(W)\Omega(U)).[s]$ and
$H\cap \Omega(U)\leq (\Omega(U_a)\Omega(U_b)).[t]$ where
$U=U_a\perp U_b$ is a decomposition into non-degenerate subspaces.
Then $H$ is a subgroup of $M_1^*=(\Omega(W\perp
U_a)\Omega(U_b)).[u]$ where $M_1^*$ is the stabilizer of the
decomposition $V=(W\perp U_a)\perp U_b$.
\end{lemma}
\begin{proof}
Again, by \cite[Lemma 4.1.1]{kl}, if $g\in M^*$ then $g=g_1g_2$
where $g_1\in I(W), g_2\in I(U)$. Furthermore, if $g\in H$ then
$g_2=g_{2a}g_{2b}$ where $g_{2a}\in I(U_a), g_{2b}\in I(U_b)$ and
$U=U_a\perp U_b$ is a decomposition of $W$ into non-degenerate
subspaces of dimension $r$ and $n-m-r$ respectively. Thus if $g\in
H$ then $g=(g_1g_{2a})g_{2b}$ and we have the required
decomposition.
\end{proof}

\begin{lemma}\label{lemma:imprim_final}
Suppose $\Lalph$ is a subgroup of maximal subgroups isomorphic to
$\hatt (X_{m_i}(q) \times Y_{n-{m_i}}(q)).[s]$ with $m_i\leq
m^\dagger\leq \frac{n}{2}$ where $m^\dagger$ is a constant. Suppose $\Lalph^*\cap \Omega(U_i)$
(where $\Omega(U_i)\cong X_{m_i}(q)$) is not irreducible in
$\Omega(U)$ then one of the following situations holds:
\begin{enumerate}
\item   There exists $j$ such that $\Lalph^* \geq \Omega(U_j)$ and $\Lalph^*$ preserves a non-degenerate decomposition, $V=W_j\perp U_j$, where
$W_j$ has dimension at most $m^\dagger$;
\item   $(V,\kappa)$ is orthogonal, $n$ is odd, $\Lalph^*\leq \Omega ^+_{n-1}(q).2$, $\Lalph^*\cap \Omega^+_{n-1}(q)$ is a
subgroup of a parabolic $P_{\frac{n-1}2}^*$ and $q\frac{n-1}{2}$
is odd;
\item   $(V,\kappa)$ is orthogonal, $n$ is even, $Y_{n-m}=\Omega ^\epsilon_{n-m}$, $\Lalph\cap
\Omega(U)\leq Sp_{n-m-2}(q)$ and $q$ is even;
\end{enumerate}
\end{lemma}
\begin{proof}
Suppose that the second and third possibilities do not occur. By
Lemma \ref{lemma:imprim_parabolics}, $\Lalph^*\cap \Omega(U)$ does
not lie in a parabolic subgroup of $\Omega(U)$. Thus $\Lalph^*\cap
\Omega(U)$ preserves a non-degenerate subspace of $U$ and so, by
Lemma \ref{lemma:imprim_nondeg}, either $\Lalph^* \geq \Omega(U)$
(and the first situation holds) or $\Lalph^* \leq I(W_0)I(U_0)$
where $V=W_0\perp U_0$ and $W_0$ has dimension $m+r$ for some
$r>0$.

We repeat our analysis using $U_0$ instead of $U$. Since $m+r\leq
m^\dagger$ this process must eventually terminate with $\Lalph^* \geq
\Omega(U_j)$ as required.
\end{proof}

Note that the first possibility of Lemma \ref{lemma:imprim_final}
has already been analysed in \stepfive. Thus, to complete our
analysis of the situation where $\Lalph$ preserves a
non-degenerate subspace of $W$, we need only consider the
irreducible subgroups of $\Omega(U)$ and the exceptions listed.

\end{itemize}

\section{The First Steps}

\subsection{Steps 1 and 2}

Steps 1 and 2 and of our method can be completed very easily. We refer to Table \ref{table:hvalues} and Table \ref{table:upperbounds} which give information about our elements $g$ and $h$ for different socles $L$ with $n\geq N_1$. We choose our elements $g$ and $h$ as follows:
\begin{itemize}
\item Suppose $nq$ is odd. When $L=PSL_n(q)$, we choose $g$ and $h$ to be involutions such that $g^*$ and $h^*$ act as the identity on a 1-dimensional non-degenerate subspace $X$ while taking all vectors in $Y$ to their negative, where $Y$ is some subspace such that $V=X\oplus Y$. When $L\neq PSL_n(q)$ we make the same choice but we require that $X$ is non-degenerate and $Y=X^\perp$.
\item When $q$ is odd and $n$ is even, we choose $g$ to be as above but this time $X$ is 2-dimensional.  If $L$ is not orthogonal then $h$ is taken to be a transvection. If $L$ is orthogonal then no transvections exist and $h$ is taken to be of the same type as $g$.
\item When $q$ is even, we take $g$ and $h$ to be unipotent elements (transvections in the non-orthogonal case and root elements in the orthogonal case).

\end{itemize}

Note that, for later steps, we may change our choice of $g$ and $h$. Now Table \ref{table:upperbounds} gives values, when $n\geq N_1$, for:
\begin{itemize}
\item $v_{as}=2w^2c$: This is an upper bound for $v$ obtained using Lemma \ref{lemma:main}. This bound applies when $G$ is almost simple.
\item $v_{s}=2w^2c$: Again this is an upper bound for $v$ obtained using Lemma \ref{lemma:main}. However this bound applies when $G$ is simple (thus, from Step 3 (Irreducibles) onwards.
\end{itemize}

In order to complete Steps 1 and 2 it is enough to observe that $|L:\hatt G(q^\frac15)|>v_{as}$ for $n\geq N_1$. Recall that $G(q^\frac15)$ is the group of similarities, of the same classical geometry as $L$, over a field of size $q^\frac15$. In certain borderline cases ($n=N_1$, say) it may be necessary to refine our value for $v_{as}$ using precise values for $w$ and $c$.

Thus, from now on, we take $G=L, n\geq N_2$ and observe that $v<v_s$ in all cases.

\begin{table}
\begin{center}
\begin{tabular}{|c|c|c|}
\hline \hline
$L$ & Conditions & $|L:C_L(h)|$ \\
\hline \hline
$PSU_n(q)$ & $n$ even, $q$ odd & $\frac{(q^n-1)(q^{n-1}+1)}{q+1}$ \\
\hline
$PSL_n(q)$ & $n$ even, $q$ odd &  $\frac{(q^n-1)(q^{n-1}-1)}{q-1}$ \\
\hline\hline
\end{tabular}
\caption{Values for $c$ when $h$ is not conjugate to $g$, $n\geq N_1$}\label{table:hvalues}
\end{center}
\end{table}

\begin{table}
\begin{center}
\begin{tabular}{|c|c|c|c|c|}
\hline \hline
$L$ & Conditions & $|L:C_L(g)|$ & $\log_q(v_{as})$ & $\log_q(v_s)$\\
\hline \hline
$\Omega_n(q)$ & $n$ odd, $q$ odd & $\frac12q^{\frac{n-1}2}(q^{\frac{n-1}2}+\epsilon)$ & 3n-3 &3n-3\\
\hline
$P\Omega^\epsilon_n(q)$ & $n$ even, $q$ odd & $\frac{q^{n-2}(q^{\frac{n-2}2}+\zeta)(q^{\frac{n}2}-\epsilon)}{2(q-\eta)}$ & 6n-5& 6n-11\\
\hline
$\Omega^\epsilon_n(q)$ & $n$ even, $q$ even & $\frac{(q^{n-2}-1)(q^{\frac{n}2}-\epsilon)(q^{\frac{n-4}2}+\epsilon)}{q^2-1}$ & 6n-5& 6n-14\\
\hline
$PSU_n(q)$ & $nq$ odd & $\frac{q^{n-1}(q^n+1)}{q+1}$ & 6n-5 & 6n-5 \\
\hline
$PSU_n(q)$ & $n$ even, $q$ odd & $\frac{q^{2n-4}(q^n-1)(q^{n-1}+1)}{(q^2-1)(q+1)}$ & 10n-5& 10n-6\\
\hline
$PSU_n(q)$ & $q$ even & $\frac{(q^n-(-1)^n)(q^{n-1}-(-1)^{n-1})}{q+1}$ & 6n& 6n-5\\
\hline
$PSp_n(q)$ & $n$ even, $q$ odd & $\frac{q^{n-2}(q^n-1)}{q^2-1}$ & 5n-1 & 5n-4\\
\hline
$Sp_n(q)$ & $n$ even, $q$ even &  $q^n-1$ & 3n+1 & 3n+1\\
\hline
$PSL_n(q)$ & $nq$ odd & $\frac{q^{n-1}(q^n-1)}{q-1}$ & 6n-2 & 6n-4\\
\hline
$PSL_n(q)$ & $n$ even, $q$ odd & $\frac{q^{2n-4}(q^n-1)(q^{n-1}-1)}{(q^2-1)(q-1)}$ & 10n-13& 10n-14\\
\hline
$PSL_n(q)$ & $q$ even & $\frac{(q^n-1)(q^{n-1}-1)}{q-1}$ & 6n+1 & 6n-2 \\
\hline
\end{tabular}
\caption{Upper bounds for $v$, $n\geq N_1$}\label{table:upperbounds}
\end{center}
\end{table}

\subsection{Steps 3 and 4}

For Steps 3 and 4 we refer to Table \ref{table:steps34}. We define
\begin{itemize}
\item $M_{Irr}$ to be the largest irreducible subgroup in $L$. Then $|L:M_{Irr}|$ is a lower bound for the index of a proper irreducible subgroup in $L$; values are obtained using \cite{liebeck3} and are valid for $n\geq\frac{N_3}{2}$. There is one exception to this (see * in the table): When $L=Sp_n(q)$, for $q$ even, there is an irreducible subgroup $O^\epsilon_n(q)$ which is the largest proper irreducible subgroup of $Sp_n(q)$. The value given at * in Table \ref{table:steps34} is a lower bound for the index of all other proper irreducible subgroups in $L$.
\item $M_{Red}$ is a reducible subgroup of $L$. When $L=PSL_n(q)$, $M_{Red}=P_m$, a parabolic subgroup preserving an $m$-dimensional subspace where $m\leq\frac{n}2$. When $L\neq PSL_n(q)$, $M_{Red}$ is of type $X_m\perp Y_{n-m}$, and so preserves a decomposition into non-degenerate subspaces. Again this value is valid for $n\geq\frac{N_3}2$.
\item $m^\dagger$ is the largest value of $m$ such that $|L:M_{Red}|\leq v_{s}$. This calculation holds for $n\geq \max\{N_1, \frac{N_3}2\}$.
\end{itemize}

In almost all cases to complete Step 3 (Irreducibles) we simply observe that $|L:M_{Irr}|>v_s$ for $n\geq N_2$. Again, for small values of $n$ ($n=N_2$, say) it may be necessary to refine our value for $v_{as}$ and $|L:M_{Irr}|$. Then to complete Step 4 (Reducibles) we read off the value for $m^\dagger$ (again, with refinement if necessary). When $L=PSL_n(q)$, no value for $m^\dagger$ is given. Instead Step 4 (Reducibles) is subsumed into Steps 5 and 6 in this case.

There is one special situation in which extra work is needed in order to complete Step 3 (Irreducibles). This occurs when $L=Sp_n(q)$, with $q$ even, and $\Lalph\leq M<L$ for $M=O_n^\epsilon(q)\in\curlyc{8}$. Then $M$ consists of the subgroup of elements which preserve a quadratic form $\kappa$ over
our vector space $V$ such that $\kappa$ polarises to a symplectic form $f$, for which $L$ is a set of linear isometries.

We need to prove that, for $n\geq 14$ our hypothesis leads to a contradiction. Then there are two possibilities:

\begin{itemize}
\item $\Lalph=M=O_n^\epsilon(q)$: Then $v=\frac12 q^{\frac{n}2}
(q^{\frac{n}2}+\epsilon)$. When $q=2,$ $L$ acts transitively on the
conjugates of $\Lalph$ (\cite{cks}; hence linear-space action will
be flag-transitive (\cite{bdd}). If we consult the list of
flag-transitive actions (\cite{ftclass}) we find that no such
actions exist and so we assume that $q>2$.

Now $O_n^\epsilon(q)>Sp_{n-2}(q)$ and so $a\geq q^{n-2}-1$. Hence
$k\leq\frac{w}{a}\leq 2(q^2+1)$ and
$b>\frac{v^2}{k^2}>\frac{1}{32}q^{2n-4}$.

Examining intersections of conjugates of $O_n^\epsilon(q)$ in $Sp_n(q)$ we find that there exist distinct conjugates which both
contain $U_P:O_{n-2}^\epsilon(q)$ where $U_P=\mathit{O}_p(P_1)$. Here $P_1$ is a parabolic subgroup of $O_n^\epsilon(q)$ which preserves a 1-dimensional non-singular subspace of the $n$-dimensional orthogonal space.
Then $b$ divides
\begin{eqnarray*}
&&(v(v-1), |Sp_n(q): (U_P:O_{n-2}(q))|)\\
&=&\frac12 q^{\frac{n}2}
(q^{\frac{n}2}+\epsilon)(q^{\frac{n}2}-\epsilon)(\frac12
(q^{\frac{n}2}+2\epsilon), q^{\frac{n-2}2}+\epsilon).
\end{eqnarray*}
For $q>2$, $(\frac12 (q^{\frac{n}2}+2\epsilon),
q^{\frac{n-2}2}+\epsilon)<\frac12 q$ and $b<\frac14
q^{\frac{3n}2+1}$. But then $\frac{1}{32}q^{2n-4}<\frac14
q^{\frac{3n}2+1}$ implies that $n<14$ as required.

\item   $\Lalph<M$: If $\Lalph$ lies in a reducible subgroup of
$M$ then $\Lalph$ lies in a reducible subgroup of $L$ and this is
dealt with in Section \ref{S:spqeven}. If $\Lalph$ lies in an irreducible subgroup of
$M$ then we can apply bounds given by \cite[Theorems 5.4 and
5.5]{liebeck3} to conclude that $|M:\Lalph|>q^{\frac14
(n^2-2n-16)}$. Since $|L:M|>q^{n-2}$ we know that
$|L:\Lalph|>{\frac14 (n^2+2n-24)}$. Since $v<q^{3n+1}$ we find
that $n\leq 12$ as required.

\end{itemize}

\begin{table}
\begin{center}
\begin{tabular}{|c|c|c|c|c|c|}
\hline \hline
$L$ & Conditions & $\log_q|L:M_{Irr}|$ & $\log_q|L:M_{Red}|$ & $m^\dagger$\\
\hline \hline
$\Omega_n(q)$ & $n$ odd, $q$ odd & $\frac14n^2-\frac14n-1$ & $mn-m^2-n+m-2$ & 3\\
\hline
$P\Omega^\epsilon_n(q)$ & $n$ even, $q$ odd & $\frac14n^2-\frac12n-\frac52$ & $mn-m^2-n-m-3$ & 14\\
\hline
$\Omega^\epsilon_n(q)$ & $n$ even, $q$ even & $\frac14n^2-\frac12n-4$ & $mn-m^2-n-m-3$ & 14\\
\hline
$PSU_n(q)$ & $nq$ odd & $\frac12 n^2+\frac12n-3$ & $2mn-2m^2-n+m$ & 4\\
\hline
$PSU_n(q)$ & $n$ even, $q$ odd & $\frac12 n^2-\frac12n-3$ & $2mn-2m^2-n+m$ & 7\\
\hline
$PSU_n(q)$ & $q$ even & $\frac13 n^2$ & $2mn-2m^2-n+m$ & 5\\
\hline
$PSp_n(q)$ & $n$ even, $q$ odd & $\frac14n^2-\frac14n-1$ & $mn-m^2$ & 6\\
\hline
$Sp_n(q)$ & $n$ even, $q$ even & * $\frac14n^2-\frac14n-1$ & $mn-m^2$ & 6\\
\hline
$PSL_n(q)$ & - & $\frac12n^2-n$ & $mn-m^2$ & - \\
\hline
\end{tabular}
\caption{Index lower bounds for $n\geq\frac{N_3}2$}\label{table:steps34}
\end{center}
\end{table}

We now proceed with Steps 5 and 6 for $L$ in various families of simple classical groups.

\section{$L=\Omega_n(q)$, $nq$ odd}

In this section we set $(N_1, N_2, N_3)=(7,13,21)$. We will prove
that $n\geq 13$ leads to a contradiction for $G$ point-primitive,
while $n\geq 21$ leads to a contradiction in all cases.

\subsection{\stepfive}



Here $M=(\Omega(U)\Omega(W)).[s]$. We have three
cases here corresponding to $m=1,2$ or $3$. Assume that
$\Omega(U)\leq \Lalph.$ By Lemma \ref{lemma:noncentral} we know
that $a\geq P(\Omega(U))$.

If $m=3$ then $M\cong(\Omega^\zeta_{n-3}(q)\times\Omega_3(q)).[4]$
and $|L:M|>q^{3n-9}$. Now $a\geq P(\Omega^\zeta_{n-3}(q))>q^{n-5}$
(\cite[Table 5.2.A]{kl}) and so $v<q^{n+7}$ which is a
contradiction for $n\geq 13$.

If $m=2$ then $M\cong(\Omega_{n-2}(q)\times\Omega^\zeta_2(q)).[4]$
and $|L:M|>q^{2n-5}.$ Now $a\geq P(\Omega_{n-2}(q))>q^{n-4}$
(\cite[Table 5.2.A]{kl}) and so $v<\frac{q^{3n-3}}{a^2}\leq
q^{n+5}$ which is a contradiction for $n\geq 13$.

If $m=1$ then $M\cong\Omega^\zeta_{n-1}(q).2$.
Then $a\geq
P(\Omega^\zeta_{n-1}(q))>q^{n-3}+q^{\frac{n-5}{2}}$ (\cite[Table
5.2.A]{kl}); this implies that $k<\frac{2w}{a}<q^2+1$ and so
$b>\frac{v^2}{k^2}>q^{2n-8}$.

Clearly there exists a distinct conjugate of $\Lalph$, $L_\beta$
say, such that $\Lalph\cap L_\beta\geq \Omega_{n-2}(q)$. Thus
$\Omega_{n-2}(q)$ fixes the line between $\alpha$ and $\beta$.
Hence $b$ divides
$|\Omega_n(q):\Omega_{n-2}(q)|=q^{n-2}(q^{n-1}-1)$.

Since $b=\frac{v(v-1)}{k(k-1)}$ divides $q^{n-2}(q^{n-1}-1)$ and
$|v|_p=q^{\frac{n-1}{2}}$ we deduce that $b$ divides
$q^{\frac{n-1}{2}}(q^{n-1}-1)$. This is a contradiction for $n\geq
13$.


\subsection{\stepsix}

Here we consider the imprimitive situation, $\Lalph<M$. Since
$k<q^{\frac{n-1}{2}}(q^{\frac{n-1}{2}}+1)$ we can apply the bounds
of \cite{deldoy} to get $|L:M|<q^{2n-2}$.

Since $|L:M|<q^{2n-2}$, $M$ must stabilize a non-degenerate
subspace $W$ of dimension at most $2$.

Suppose that $W$ has dimension $2$. Then $M\cong
(\Omega_2^\zeta(q)\times \Omega_{n-2}(q)).[4]$, $|L:M|>q^{2n-5}$ and
$|M:\Lalph|<q^{n+2}$. By Lemma \ref{lemma:imprim_final}
$\Lalph\cap\Omega(U)$ is irreducible. But, by \cite[Theorem
5.6]{liebeck3}, this implies that $|M:\Lalph|>q^{n+2}$.

Suppose that $W$ has dimension $1$. Then $M\cong
\Omega_{n-1}^\zeta(q).2$,
$|L:M|=\frac{1}{2}q^{\frac{n-1}{2}}(q^{\frac{n-1}{2}}+\zeta)$ and
$|M:\Lalph|<q^{2n-1}$. By Lemma \ref{lemma:imprim_final} there are
two possibilities:
\begin{itemize}
\item $\Lalph\cap\Omega(U)$ is irreducible. But, by \cite[Theorem
5.6]{liebeck3}, this implies that $|M:\Lalph|>q^{2n-1}$ for $n\geq
21$.
\item $\epsilon=+$ and $\Lalph\cap\Omega(U)$ lies in a parabolic
subgroup $P_{\frac{n-1}{2}}$ (see Lemma
\ref{lemma:imprim_final}). But then
$|M:\Lalph|>q^{\frac18(n-1)(n-3)}>q^{2n-1}$ for $n\geq 21$.
\end{itemize}

\section{$L=P\Omega^\epsilon_n(q)$, $n$ even, $q$ odd}\label{section:omega n even q odd}

In this section we set $(N_1,N_2,N_3)=(18,26,26)$. We will prove
that $n\geq 26$ leads to a contradiction.

\subsection{\stepfive}

Suppose first of all that $\Lalph$ lies inside a maximal subgroup
of type $O_m\perp O_{n-m}$ with $m\leq 14$. Then, by assumption, $\Lalph>\hatt \Omega^\zeta_{n-14}(q)$ and $\Lalph$ contains at least
$$\frac{q^{n-16}(q^{\frac{n-14}{2}}-1)(q^{\frac{n-16}{2}}-1)}{4(q+1)}$$
conjugates of $g$. Then, by Lemma \ref{lemma:main}, we have $k<
\frac{2w}a\leq 2q^{31}$ and $v<q^{2n+59}$. This implies that
$$mn-n-m^2-m-3<2n+59.$$
Since $2m\leq n$ and $n\geq 26$ we conclude that $m\leq 9$. But this implies that $\Lalph>\Omega^\zeta_{n-10}(q)$ and we repeat the process to find that $v<q^{2n+43}$. But this implies that $\Lalph>\Omega^\zeta_{n-8}(q)$, $v<q^{2n+35}$ and $m\leq 6$. We continue to find that $v<q^{2n+31}$ and $m\leq 2$.

When $m=2$ we find that $q^{2n-6}<v$ and $|v|_p=q^{n-2}$. Clearly there exists a distinct
conjugate of $\Lalph$, $L_\beta$ say, such that $\Lalph\cap
L_\beta\geq \Omega^\eta_{n-4}(q)$. This must lie in a
line-stabilizer hence we conclude that
$$b<q^{n-2}|\Omega^\epsilon_n(q): \Omega^\eta_{n-4}(q)|_{p'}<q^{3n-5}.$$
Now, since $\Lalph>\Omega^\zeta_{n-2}(q)$, $k<\frac{2w}{a}<2q^3$.
Hence $b>\frac{v^2}{k^2}>q^{4n-20}$. But this gives a
contradiction for $n\geq 26$.

When $m=1$, $\Lalph\cong \Omega_{n-1}(q)$ and
$v=\frac{1}{2}q^{\frac{n-2}{2}}(q^{\frac{n}2}-1)>q^{n-3}$. Then,
as before, $b>\frac{v^2}{k^2}>\frac14 q^{2n-12}$. Now, similarly
to before, we know that a line-stabilizer must contain
$\Omega^\epsilon_{n-2}(q)$ and so
$$b<q^{\frac{n-2}2}|\Omega^\epsilon_n(q): \Omega^\epsilon_{n-2}(q)|\leq q^{\frac{n-2}{2}}(q^{\frac{n}2}+1)(q^{\frac{n-2}2}+1).$$
Once again this yields a contradiction for $n\geq 26$.

\subsubsection{Other primitive possibilities}

Suppose that $\epsilon=+$ and $\Lalph$ stabilizes a totally
singular subspace of dimension $\frac{n}2$. So
$\Lalph\cong\hatt[q^{\frac18n(n-2)}]:\frac12 GL_{\frac{n}2}(q)$. Then
$$a\geq \frac{q^{\frac18n(n-2)}\times\frac12 \times |GL_{\frac{n}2}(q)|}{q^{\frac18(n-2)(n-4)}\times\frac12\times |GL_{\frac{n-2}2}(q)|\times 4 \times |\Omega^-_2(q)|}.$$
This implies that $k\leq \frac{2w}{a}<q^{\frac{n}2+1}$ and
$v<q^{3n}$. But $v>q^{\frac18n(n-2)}$ and we have a
contradiction for $n\geq 26$.

\subsection{\stepsix}

Here we consider the imprimitive situation, $\Lalph<M$. We know that $k<2q^{2n-4}$ and we can
apply the bounds of \cite{deldoy} to get $|L:M|<\frac12k^2< q^{4n-7}$.

Suppose that $M$ is a maximal subgroup of type $O_m\perp O_{n-m}$. This implies that
$$mn-n-m^2-m-3<4n-7.$$
Since $2m\leq n$ and $n\geq 26$, we conclude that $m\leq 6$. By Lemma \ref{lemma:imprim_final}, $\Lalph^*\cap\Omega(U)$ is irreducible. Then Theorems 5.4 and 5.5 of \cite{liebeck3} imply that $|\Omega(U):\Lalph^*\cap\Omega(U)|>q^{\frac14((n-m)^2-2(n-m)-10)}$ and $v>q^{6n-11}$ for $n\geq 32$. Checking the indices of maximal irreducible subgroups for $n=26,28$ and $30$ we are able to exclude these also.

\subsubsection{Other imprimitive possibilities}

Suppose that $\epsilon=+$. If $\Lalph$ lies in a parabolic subgroup $P_{\frac{n}2}$ then
$$\frac18n(n-2)<4n-7.$$
This implies that $n\leq 32$. We need only consider the situation
where $\frac{n}2$ is odd hence we are left with $n=30$ and $n=26$.

Let $U={\it O}_p(P_{\frac{n}2})$. We consider two situations:
\begin{itemize}
\item {\bf Suppose that $\Lalph\cap U<U$}. Then we can apply Proposition \ref{prop:deldoy} to $U:(\Lalph\cap \hatt SL_{\frac{n}2}(q))$. Thus this subgroup must have index at most $q^{4n-7}$. When $n=30$ this implies that $\Lalph\geq \hatt SL_{15}(q)$. But then $a>q^{15}$ and $\frac18n(n-2)>4n-37$ which is a contradiction. When $n=26$ we know that $\Lalph\cap \hatt SL_{13}(q))$ must lie inside a parabolic subgroup, $P_1$. Iterating the procedure we find that $\Lalph\geq SL_{12}(q)$ and $a>q^{12}$ which yields a contradiction.

\item {\bf Suppose that $\Lalph\geq U$}. Then let $g$ and $h$ be root elements and observe that $U$ contains more than $q^{2n}$ root elements. Then $w=c=\frac{(q^{n-2}-1)(q^{\frac{n}2}-1)(q^{\frac{n-4}2}+1)}{q^2-1}$ and, by Lemma \ref{lemma:main}, $v<q^{2n-10}$. This gives a contradiction for $n=26$ and $n=30$.

\end{itemize}

\section{$L=\Omega^\epsilon_n(q)$, $n$ even, $q$ even}\label{section:omega n even q even}

In this section we set $(N_1,N_2,N_3)=(16,26,26)$. We will prove
that $n\geq 26$ leads to a contradiction.

\subsection{\stepfive}

Suppose first of all that $\Lalph$ is a subgroup of a maximal
subgroup of type $O_m\perp O_{n-m}$ with $m$ even. Now, by assumption, $\Lalph>\Omega^\zeta_{n-m}(q)$. Then $a>q^{2n-2m-7}$ and so $v<q^{2n+4m}$. This implies that
$$mn-n-m^2-m-3<2n+4m.$$
Since $2m\leq n$ and $n\geq 26$ we conclude that $m\leq 4$.

When $m=4$ we find that $v<q^{2n+16}$ and, examining $|L:\Lalph|,$
we obtain a contradiction for $n\geq 26$. When $m=2$ we find that
$q^{2n-6}<v$ and $|v|_p\leq q^{n-2}$. As in the odd characteristic
case $b<q^{3n-5}.$ Since $\Lalph>\Omega^\zeta_{n-2}(q)$ we have
$k<\frac{2w}{a}<q^8$ and $b>\frac{v^2}{k^2}>q^{4n-28}$. But this
gives a contradiction for $n\geq 26$.

\subsubsection{Other primitive possibilities}

Suppose that $\epsilon=+$ and $\Lalph$ stabilizes a totally
singular subspace of dimension $\frac{n}2$. So
$\Lalph\cong[q^{\frac18n(n-2)}]: GL_{\frac{n}2}(q)$ and it
is easy to see that $a\geq q^{\frac{3n}{2}-6}$. By Lemma \ref{lemma:main}, we have that $v<q^{3n-2}$. But
$v>q^{\frac18n(n-2)}$ and we have a contradiction for $n\geq 26$.

Finally consider the possibility that $\Lalph=Sp_{n-2}(q)$ and
$v=q^{\frac{n-2}2}(q^{\frac{n}2}-\epsilon)$. Now $\Lalph$ is the
set of elements in $L$ which stabilizes a $1$-dimensional
non-singular subspace of $V$. Consider
$H=\Omega_2^\zeta(q))\times\Omega_{n-2}^\epsilon(q)$ stabilizing a
$2$-dimensional non-degenerate subspace of $V$. Then it is easy to
see that $\Omega_{n-2}^\epsilon(q),$ normal in $H$, stabilizes 2
distinct $1$-dimensional non-singular subspaces of $V$. We
conclude that $\Omega_{n-2}^\epsilon(q)$ stabilizes a line and so
$b$ divides
$$|v|_p|\Omega^\epsilon_n(q):\Omega^\epsilon_{n-2}(q)|_{p'}=q^{\frac{n-2}2}(q^{\frac{n}2}-\epsilon)(q^{\frac{n-2}2}+\epsilon)<q^{\frac{3n}2-1}.$$

Now observe that, since $\Omega_{n-2}^\epsilon(q)$ stabilizes a point,
$$w\leq\frac{(q^{n-2}-1)(q^{\frac{n}2}-\epsilon)(q^{\frac{n-4}2}+\epsilon)}{q^2-1}, \ \ a\geq \frac{(q^{n-4}-1)(q^{\frac{n-2}2}-\epsilon)(q^{\frac{n-6}2}+\epsilon)}{q^2-1}.$$
Applying our bound we find that
$$k<\frac{2w}a\leq \frac{2(q^{\frac{n-2}2}+\epsilon)(q^{\frac{n}2}-\epsilon)}{(q^{\frac{n-4}2}-\epsilon)(q^{\frac{n-6}2}+\epsilon)}\leq q^6$$
and so $b>\frac{v^2}{k^2}>q^{2n-13}$. But then $2n-13<\frac{3n}2-1$ implies that $n<26$ as required.

\subsection{\stepsix}

Here we consider the imprimitive situation,
$\Lalph<M$. We know that $k<2q^{2n-4}$ and so we can apply the bounds of \cite{deldoy} to get $|L:M|<2w^2<q^{4n-9}$.

Suppose that $M$ is a maximal subgroup of type $O_m\perp O_{n-m}$
with $m$ even. This implies that
$$mn-n-m^2-m-3<4n-9.$$
Since $2m\leq n$ and $n\geq 26$, we conclude that $m\leq 6$. In
what follows we take $m$ to be the largest integer such that
$\Lalph$ is contained in a maximal subgroup of type $O_m\perp
O_{n-m}$; thus $m=2,4$ or $6$. We refer to Lemma
\ref{lemma:imprim_final} and go through all possibilities:

\begin{itemize}
\item   $\Lalph\cap \Omega(U)$ is irreducible. In this case
$|M:\Lalph|\geq q^{\frac14((n-m)^2-2(n-m)-16)}$. But then
$|L:\Lalph|>q^{6n-14}$ for $n\geq 26$.

\item $\Lalph\cap \Omega(U)\leq Sp_{n-m-2}(q)$. In this case
$$\Lalph\leq M_1=(\Omega(W)\times Sp_{n-m-2}(q)).2$$ and
$|L:M_1|>q^{mn+n-m^2-m-4}$. Once again we have a number of
possibilities:
\begin{itemize}
\item  $\Lalph\cap Sp_{n-m-2}(q)$ is irreducible:
\begin{itemize}
\item  $\Lalph\cap Sp_{n-m-2}(q)\geq \Omega^\epsilon_{n-m-2}(q)$. Then, similarly to before, $v<q^{2n+4(m+2)}$ which is a contradiction for $n\geq 26$.
\item  $\Lalph\cap Sp_{n-m-2}(q)\leq O^\epsilon_{n-m-2}(q)$. Then $\Lalph<(\Omega_{m+2}^{\eta}(q)\times\Omega_{n-m-2}^{\epsilon}(q)).2$ which is a contradiction.
\item  $\Lalph\cap Sp_{n-m-2}(q)$ lies in any other irreducible subgroup. Then, referring to \cite[Table 3.5.C]{kl}, we find that
$|M_1:\Lalph|>q^{\frac14(n-m-2)^2-\frac14(n-m-2)-1}$ and $|L:\Lalph|>q^{6n-14}$ for $n\geq 26$.
\end{itemize}
\item  $\Lalph\cap Sp_{n-m-2}(q)$ lies inside $P$ a parabolic
subgroup of $Sp_{n-m-2}(q)$. This possibility can be excluded by
an argument similar to the proof for Lemma
\ref{lemma:imprim_parabolics}.
\item  $\Lalph\cap Sp_{n-m-2}(q)\leq Sp_{r}(q)\times Sp_{n-m-2-r}(q)$.
Examining
$$|L:(\Omega(W)\times Sp_r(q)\times Sp_{n-m-2-r}(q)).2|$$ we
conclude that $m+r\leq 8$. In a similar way to Lemma
\ref{lemma:imprim_final} we need only examine the possibility that
$\Lalph\cap Sp_{n-m-2-r}(q)$ is irreducible. There are a number of
possibilities:
\begin{itemize}
\item  $\Lalph\geq Sp_{n-m-2-r}(q)\geq
\Omega^\epsilon_{n-m-2-r}(q)$. Then, as before, $v<q^{2n+4(m+r+2)}$ which is a contradiction for $n\geq 26$.
\item
$\Lalph\cap Sp_{n-m-2-r}(q)\leq O_{n-m-2-r}(q)$. Then
$\Lalph<(\Omega_{m+2+r}^{\eta}(q)\times
\Omega_{n-m-2-r}^{\zeta}(q)).2$ which is a contradiction. \item
$\Lalph\cap Sp_{n-m-2-r}(q)$ lies in any other irreducible
subgroup. Then, referring to \cite[Table3.5.C]{kl}, we find that
$$|Sp_{n-m-2-r}(q):\Lalph\cap
Sp_{n-m-2-r}(q)|>q^{\frac14(n-m-r-2)^2-\frac14(n-m-2-r)-1}$$ and
$|L:\Lalph|>q^{6n-14}$ for $n\geq 26$.
\end{itemize}
\end{itemize}

\end{itemize}

\subsubsection{Other imprimitive possibilities}

Suppose first of all that $\epsilon=+$ and $\Lalph$ lies in a parabolic subgroup $P_{\frac{n}2}$. Then
$$\frac18n(n-2)<4n-9$$
and so $n\leq 32$. We need only consider the situation
where $\frac{n}2$ is odd hence we are left with $n=30$ and $n=26$. These cases are ruled out very similarly to when $q$ is odd.

Now suppose that $\Lalph<M=Sp_{n-2}(q)$. Then $|L:M|>q^{n-2}$ and so
$|M:\Lalph|<q^{5n-12}$. If $\Lalph<O_{n-2}(q)$ then $\Lalph$
preserves a non-degenerate subspace of dimension $2$ and this is
covered above. If $\Lalph$ lies in any other irreducible subgroup
of $Sp_{n-2}(q)$ then $|M:\Lalph|>q^{5n-12}$ which is a
contradiction for $n\geq 26$.

By a similar argument to Lemma \ref{lemma:parabolics} we can
conclude that $\Lalph$ does not lie in a parabolic subgroup of
$Sp_{n-2}(q)$. Thus $\Lalph$ must lie in a maximal subgroup $M_1$
of type $Sp_m\perp Sp_{n-m-2}$. Since $|M:\Lalph|<q^{5n-12}$ we conclude that $m\leq 6$.

In a similar way to Lemma \ref{lemma:imprim_final} we need only examine the
possibility that $\Lalph\cap Sp_{n-m-2}(q)$ is irreducible. There
are a number of possibilities:
\begin{itemize}
\item  $\Lalph\geq Sp_{n-m-2}(q)\geq \Omega^\epsilon_{n-m-2}(q)$. Once more this implies that $v<q^{2n+4(m+2)}$ which is a contradiction for $n\geq 26$.
\item  $\Lalph\cap Sp_{n-m-2}(q)\leq O_{n-m-2}(q)$. Then $\Lalph<(\Omega_{m+2}^{\eta}(q)\times
\Omega_{n-m-2}^{\zeta}(q)).2$ which is covered above.
\item
$\Lalph\cap Sp_{n-m-2}(q)$ lies in any other irreducible subgroup.
Then, referring to \cite[Table3.5.C]{kl}, we find that
$|Sp_{n-m-2}(q):\Lalph\cap Sp_{n-m-2}(q)|>q^{\frac14(n-m-2)^2-\frac14(n-m-2)-1}$
and $|L:\Lalph|>q^{6n-14}$ for $n\geq 26$.
\end{itemize}

\section{$L=PSU_n(q)$, $nq$ odd}

In this section we set $(N_1, N_2, N_3)=(11,15,15)$. We will prove
that $n\geq 15$ leads to a contradiction.

\subsection{\stepfive}

Here $\Lalph\leq\hatt ((SU_m(q)\times SU_{n-m}(q)).(q+1))$ for
some $m\leq 4$. Observe that in all cases $\Lalph>\hatt
SU_{n-6}(q)$ and so $a\geq q^{n-5}\frac{q^{n-4}+1}{q+1}$. Thus
$k\leq \frac{2w}a\leq q^{9}$ and so $v<q^{2n+16}$. This implies
that
$$2mn-2m^2-n+m<2n+16.$$
Since $2m<n$ and $n\geq 15$ we conclude that $m\leq 2$. But this
implies that $\Lalph>\hatt SU_{n-2}(q)$ and so $a\geq
q^{n-3}\frac{q^{n-2}+1}{q+1}$ and $k<q^5$. Repeating the process
we find that $v<q^{2n+8}$ and $m=1$.

We have $\hatt SU_{n-1}(q)\leq \Lalph\leq \hatt GU_{n-1}(q)$. Now,
by Lemma \ref{lemma:selfnormalizing}, $\Lalph=\hatt GU_{n-1}(q)$,
$k<q^5$ and $v= q^{n-1}\frac{q^n+1}{q+1}$. Then
$b>\frac{v^2}{k^2}>q^{4n-16}$. Clearly there exists a distinct
conjugate of $\Lalph$, $L_\beta$ say, such that $\Lalph\cap
L_\beta\geq \hatt SU_{n-2}(q)$. This must lie in a line-stabilizer
hence we conclude that $b\big| q^{2n-3}(q^n+1)(q^{n-1}-1)$. Since
$v=q^{n-1}\frac{q^n+1}{q+1}$ we must have $b\big|
q^{n-1}(q^n+1)(q^{n-1}-1)$. But then $q^{4n-16}<b<q^{3n-3}$ which
is a contradiction.

\subsection{\stepsix}

Here we consider the imprimitive situation, $\Lalph<M$. We know that $k<2q^{2n-1}$ and we can
apply the bounds of \cite{deldoy} to get $|L:M|<q^{4n-1}$. This implies that
$$2mn-2m^2-n+m<4n-1.$$
Once again, since $2m<n$ and $n\geq 15$, we conclude that $m\leq
2$.

If $m=2$ then $|L:M|>q^{4n-9}$ and so $|M:\Lalph|<q^{2n+4}$. We examine $\Lalph\cap\Omega(U)$; by Lemma \ref{lemma:imprim_final} this must be an irreducible subgroup of $\Omega(U)$. But then, by \cite[Theorem 5.3]{liebeck3}, $|\Omega(U):\Lalph\cap\Omega(U)|>q^{2n+4}$ which is a
contradiction.

If $m=1$ then $|L:M|>q^{2n-3}$ and so $|M:\Lalph|<q^{4n-2}$. Once again $\Lalph\cap\Omega(U)$ must be an irreducible subgroup of $\Omega(U)$. But then, by \cite[Theorem 5.3]{liebeck3}, $|\Omega(U):\Lalph\cap\Omega(U)|>q^{4n-2}$ which is a contradiction.

\section{$L=PSU_n(q)$, $n$ even, $q$ odd}

In this section we set $(N_1, N_2, N_3)=(16,22,22)$. We will prove
that $n\geq 22$ leads to a contradiction.

\subsection{\stepfive}

Here $\Lalph\leq \hatt ((SU_m(q)\times SU_{n-m}(q)).(q+1))$ for
some $m\leq 7$. Observe that in all cases $\Lalph>\hatt
SU_{n-8}(q)$ and so
$$a\geq q^{2n-20}(q^{n-10}+\dots + q^2+1)(q^{n-10}-\dots-q+1).$$
Thus $k\leq \frac{2w}a\leq
2q^{26}(q^{8}+1)$ and so $v<q^{2n+64}$. This implies that
$$2mn-2m^2-n+m<2n+64.$$
Since $2m<n$ and $n\geq 22$ we conclude that $m\leq 4$. But this
implies that $\Lalph>\hatt SU_{n-4}(q).$ We repeat the process to
find that $v<q^{2n+32}$ and $m\leq 2$. But this implies that
$\Lalph>\hatt SU_{n-2}(q).$ Again we repeat the process to find
that $v<q^{2n+16}$ and $m=1$.

Thus $\hatt SU_{n-1}(q)\leq \Lalph\leq\hatt GU_{n-1}(q)$. By Lemma
\ref{lemma:selfnormalizing} $\Lalph=\hatt GU_{n-1}(q)$, $k<q^5$
and $v=q^{n-1}\frac{q^n-1}{q+1}$. Then
$b>\frac{v^2}{k^2}>q^{4n-18}$. Clearly there exists a distinct
conjugate of $\Lalph$, $L_\beta$ say, such that $\Lalph\cap
L_\beta\geq \hatt SU_{n-2}(q)$. This must lie in a line-stabilizer
hence we conclude that $b\big| q^{2n-3}(q^n-1)(q^{n-1}+1)$. Since
$v=q^{n-1}\frac{q^n+1}{q+1}$ we must have $b\big|
q^{n-1}(q^n-1)(q^{n-1}+1)$. But then $q^{4n-18}<b<q^{3n-3}$ which
is a contradiction.

\subsection{\stepsix}

Here we consider the imprimitive situation, $\Lalph<M$. We know that $k<2q^{3n-6}(q^{n-2}+\dots
q^2+1)$ and we can apply the bounds of \cite{deldoy} to get $|L:M|<\frac12k^2< q^{8n-15}$. This implies that
$$2mn-2m^2-n+m<8n-15.$$
Since $2m<n$ and $n\geq 22$, we conclude that $m\leq 4$ or $m=5$
and $n\leq 30$.

Now $|L:M|>q^{2mn-2m^2-n-m}$ . We examine $\Lalph\cap\Omega(U)$;
by Lemma \ref{lemma:imprim_final} this must be an irreducible
subgroup of $\Omega(U)$. But then, by \cite[Theorem
5.3]{liebeck3},
$|\Omega(U):\Lalph\cap\Omega(U)|>q^{\frac12(n-m)^2-\frac12(n-m)-3}$
which is a contradiction for $n\geq 26$. A closer examination of
bounds for $n=20$ and $22$ also yields a contradiction as
required.

\section{$L=PSU_n(q)$, $q$ even}

In this section we set $(N_1, N_2, N_3)=(11,13,13)$. We will prove
that $n\geq 13$ leads to a contradiction.

\subsection{\stepfive}

Here $\Lalph\leq\hatt ((SU_m(q)\times SU_{n-m}(q)).(q+1))$ for
some $m\leq 5$. Observe that in all cases $\Lalph>\hatt
SU_{n-5}(q)$ and so $a\geq
\frac{(q^{n-5}-(-1)^{n-5})(q^{n-6}-(-1)^{n-6})}{q+1}$. Thus $k\leq
\frac{2w}a\leq 2q^5(q^5+1)$ and so $v<q^{2n+20}$. This implies
that
$$2mn-2m^2-n+m<2n+20.$$
Since $2m<n$ and $n\geq 13$ we conclude that $m\leq 2$. But this
implies that $\Lalph>\hatt SU_{n-2}(q).$ We repeat the process to
find that $v<q^{2n+7}$ and $m=1$ or $(m,n,q)=(2,13,2)$. This
last situation is easily excluded.

Thus $\hatt SU_{n-1}(q)\leq \Lalph\leq\hatt GU_{n-1}(q)$. By Lemma
\ref{lemma:selfnormalizing} $\Lalph=\hatt GU_{n-1}(q)$, $k<q^3$
and $v=q^{n-1}\frac{q^n\pm1}{q+1}$. Then
$b>\frac{v^2}{k^2}>q^{4n-12}$. Clearly there exists a distinct
conjugate of $\Lalph$, $L_\beta$ say, such that $\Lalph\cap
L_\beta\geq \hatt SU_{n-2}(q)$. This must lie in a line-stabilizer
hence we conclude that $b\big| q^{2n-3}(q^n\pm1)(q^{n-1}\mp1)$.
Since $v=q^{n-1}\frac{q^n\pm1}{q+1}$ we must have $b\big|
q^{n-1}(q^n\pm1)(q^{n-1}\mp1)$. But then $q^{4n-12}<b<q^{3n-3}$
which is a contradiction.

\subsection{\stepsix}

Here we consider the imprimitive situation, $\Lalph<M$. We know that $k<2q^{2n-2}$ and we can
apply the bounds of \cite{deldoy}to get $|L:M|<\frac12k^2< q^{4n-3}$. This implies that
$$2mn-2m^2-n+m<4n-3.$$
Since $2m<n$ and $n\geq 14$, we conclude that $m\leq 2$.

Now $|L:M|>q^{2mn-2m^2-n-m}$ . We examine $\Lalph\cap\Omega(U)$;
by Lemma \ref{lemma:imprim_final} this must be an irreducible
subgroup of $\Omega(U)$. But then, by \cite[Theorem
5.3]{liebeck3},
$|\Omega(U):\Lalph\cap\Omega(U)|>q^{\frac13(n-m)^2}$ which is a
contradiction for $n\geq 20$. A closer examination of bounds for
$13\leq n\leq 19$ also yields a contradiction as required.

\section{$L=PSp_n(q)$, $n$ even, $q$ odd}

In this section we set $(N_1, N_2, N_3)=(12,22,22)$. We will prove
that $n\geq 22$ leads to a contradiction.

\subsection{\stepfive}

Here $\Lalph\leq\hatt (Sp_m(q)\circ Sp_{n-m}(q))$ for some even
$m\leq 6$. Observe that in all cases $\Lalph>Sp_{n-6}(q)$ and so
$a\geq q^{n-8}(q^{n-8}+\dots + q^2+1)$. Thus $v<q^{n+26}$. This
implies that
$$mn-m^2<n+26.$$
Since $2m<n$ and $n\geq 22$ we conclude that $m= 2$. But this
implies that $\Lalph>Sp_{n-2}(q).$ We repeat the process to find
that $v<q^{n+10}$ which is a contradiction.

\subsection{\stepsix}

Here we consider the imprimitive situation, $\Lalph<M$. We know that $k<2q^{n-2}(q^{n-2}+\dots q^2+1)$ and we can
apply the bounds of \cite{deldoy} to get $|L:M|<\frac12k^2< q^{4n-7}$. This implies that
$$mn-m^2<4n-7.$$
Since $2m<n$ and $n\geq 22$, we conclude that $m\leq 4$.

Now $|L:M|>q^{mn-m^2}$. We examine $\Lalph\cap\Omega(U)$; by Lemma
\ref{lemma:imprim_final} this must be an irreducible subgroup of
$\Omega(U)$. But then, by \cite[Theorem 5.3]{liebeck3},
$|\Omega(U):\Lalph\cap\Omega(U)|>q^{\frac14(n-m)^2-\frac14(n-m)-1}$
which is a contradiction for $n\geq 22$.

\section{$L=Sp_n(q)$, $n$ even, $q$ even}\label{S:spqeven}

In this section we set $(N_1, N_2,N_3)=(8,14,14)$. We will prove
that $n\geq 14$ leads to a contradiction.

\subsection{\stepfive}

Here $\Lalph\leq\hatt (Sp_m(q)\circ Sp_{n-m}(q))$ for some even
$m\leq 6$. Observe that in all cases $\Lalph>Sp_{n-6}(q)$ and so
$a\geq q^{n-6}-1$. Thus $v<q^{n+16}$. This implies that
$$mn-m^2<n+16.$$
Since $2m<n$ and $n\geq 14$ we conclude that $m=2$. But this
implies that $\Lalph>Sp_{n-2}(q).$ We repeat the process to find
that $v<q^{n+10}$ which is a contradiction.

\subsection{\stepsix}

Here we consider the imprimitive situation, $\Lalph<M$. We know that $k<2q^{n}$ and we can
apply the bounds of \cite{deldoy} to get $|L:M|<\frac12k^2< q^{2n+1}$. This implies that
$$mn-m^2<2n+1.$$
Since $2m<n$ and $n\geq 14$, we conclude that $m=2$.

Now $|L:M|>q^{2n-4}$. Examining \cite[Table 5.2.A]{kl} we find
that $P(Sp_{n-2}(q))>q^5$ in all cases. Thus $|L:(Sp_2(q)\times M_1)|>q^{2n+1}$ where $M_1$ is maximal in $Sp_{n-2}(q)$. By Proposition \ref{prop:deldoy}, this implies that $\Lalph\cap\Omega(U)=M_1$ and so is and so is maximal in $\Omega(U)$.

Observe now that $|M:\Lalph|<q^{n+5}$. We examine $\Lalph\cap\Omega(U)$; by Lemma \ref{lemma:imprim_final} this must
be an irreducible subgroup of $\Omega(U)$. There are two possibilities:
\begin{itemize}
\item $\Lalph=O_{n-2}^\epsilon(q)>Sp_{n-4}(q)$. Then $a\geq
q^{n-4}-1$ and $v<q^{n+11}$ which is a contradiction. \item
$\Lalph\cap\Omega(U)$ lies in any other irreducible subgroup of
$Sp_{n-2}(q)$. But in this case $|L:\Lalph|>q^{3n+1}$ which is a
contradiction.

\end{itemize}

\section{$L=PSL_n(q)$}
In this section $\Lalph$ lies in a parabolic subgroup of $L$. We will prove that, if $n\geq N_2$, then only Example \ref{example:flagtransitive} occurs.

\subsection{\stepfive}
Suppose that $H\leq\Lalph\leq P_m$ where $H\cong \hatt
SL_{n-m}(q)$ and $H$ is normal in a Levi-complement of $P_m$. Clearly, in all cases, $\Lalph$ contains transvections. Thus we take $g$ and $h$ to be transvections and observe that $a\geq (q^{n-m-1}-1)(q^{n-m-1}+\dots+q+1)$. By Lemma \ref{lemma:main} we have $v<q^{2n+4m+2}$. Since $v>q^{m(n-m)}$ we conclude that $m\leq 3$.

If $\Lalph<P_m$ then we can apply the bounds of Delandtsheer and
Doyen. Thus $|L:P_m|<\frac{2w^2}{a^2}<2(q^3+1)^4$. This is a contradiction
for $n\geq 17$.

If $\Lalph=P_1$ then this action is $2$-transitive, hence
flag-transitive (\cite{bdd}) and so corresponds to the known
action on $PG(n-1,q)$.

If $\Lalph=P_m$, for $m=2$ or $3$, then there exists a conjugate $P_m^x$ not equal to $P_m$ such that $\hatt SL_{n-m-1}(q)<P_m\cap P_m^x$. Hence $b$ divides
$|SL_n(q):SL_{n-m-1}(q)|$. We know that $p$ is not significant
(\cite {gill3}) and $p$ does not divide $v$. Hence $b$ divides
$|SL_n(q):SL_{n-m-1}(q)|_{p'}<q^{(m+1)(n-\frac12m)}$. Now $v>q^{m(n-m)}$
and $k<q^{2m+1}$ hence $b>\frac{v^2}{k^2}>q^{2mn-2m^2-4m-4}$. This
gives a contradiction for $n\geq 17$.

\subsection{\stepsix}
We are left with the possibility that $\Lalph<P_m$ and $\Lalph$
does not contain $H$. Our upper bound for $v$ varies in this case. The worst case scenario is when $n$ is even and $q$ is odd in which case $v<q^{10n-14}$ and, by Proposition \ref{prop:deldoy}, $|L:P_m|<q^{8n-14}$.

Let $U=\mathit{O}_p(P_m)$. We have two situations:

\begin{itemize}
\item {\bf Suppose that $\Lalph$ does not contain $U$}. Then, by Proposition \ref{prop:deldoy} we can apply the Delandtsheer-Doyen
bound to the group $M_1 = U\Lalph$. Then we have $|L:P_m|\cdot |H:\Lalph\cap H|<q^{8n-14}.$ We
suppose that $\Lalph\cap H$ is an irreducible subgroup of $H$.
Referring to \cite[Theorem 5.1]{liebeck3} we see that this implies
that $$|H:\Lalph\cap H|\geq q^{\frac12 (n-m)^2-(n-m)}.$$ But then
the Delandtsheer-Doyen bound is violated for $n\geq 20$.

For $q$ even or $nq$ odd the Delandtsheer-Doyen bound implies that
$|L:P_m|\cdot |H:\Lalph\cap H|<q^{4n}.$ This implies that, for $n\geq 17$, $\Lalph\cap H$ is not an irreducible subgroup of $H$.

Thus $\Lalph\cap H$ lies in a parabolic subgroup of $H$. We can continue to iterate in this way by considering $\Lalph\cap\hatt SL_{n-s}(q)$ for some $s$. This will either terminate with $\Lalph\geq \hatt SL_{n-s}(q)$ for some $s$ or it will increase the cumulative index by at least $q^s$ and increase $s$ by at least $1$. It is easy to see that this must terminate with a final value for $s\leq n-2$.

Thus $\Lalph>\hatt SL_2(q)$ and so $\Lalph$ must contain transvections and, once more, we have $|L:M_1|<q^{4n}$. Repeating our analysis with this new bound we find that $s\leq 10$. In fact, for $\Lalph>\hatt SL_{n-s}(q)$, we have $|L:M_1|<q^{4s+4}$. This implies that $s\geq n-4$ which is a contradiction for $n\geq 17$.

\item{\bf Suppose that $\Lalph\geq U$}. Then $\Lalph$ contains transvections and we can take
$w=c=(q^{n-1}-1)(q^{n-1}+\dots+q+1)$. Then $a>q^{n-1}-1$ and so
$v<q^{4n}$. This implies that $m\leq 6$. Once again we can exclude the possibility that $\Lalph\cap H$ is an irreducible subgroup of $H$.

Thus $H\cap\Lalph\leq P_{m_1}$, a parabolic subgroup of $H$. Since $v<q^{4n}$ we must have $m+m_1\leq 6$. If $H\cap\Lalph=P_{m_1}$ then observe that $\Lalph\geq \hatt SL_{n-6}(q)$.

$H\cap\Lalph<P_{m_1}$ then we can apply Proposition \ref{prop:deldoy} to $U:P_{m_1}$. This implies that $q^{n(m+m_1)-(m+m_1)^2}<|L:(U:P_{m_1})|<q^{2n+1}$. Then $m+m_1=2$ and $m=m_1=1$. We can iterate this procedure once again by examining $H_1\cap \Lalph$ where $H_1$ is normal in a Levi-complement of $P_{m_1}$. Clearly $H_1\cap\Lalph$ must lie in a parabolic subgroup of $H_1$. In fact, by considering Proposition \ref{prop:deldoy}, it must equal such a parabolic subgroup and once more $\Lalph\geq \hatt SL_{n-6}(q)$.

Now $a \geq (q^{n-7}-1)(q^{n-7}+\dots+q+1)$ and $v<q^{2n+28}$. Then $\Lalph^*$ must stabilize a subspace of dimension at most $4$ and $\Lalph>\hatt SL_{n-4}(q)$. Repeating the process we conclude that $\Lalph>\hatt SL_{n-3}(q)$ and $|L:M|<2(q^3+1)^4$ which is a contradiction for $n\geq 17$.

\end{itemize}

\section{Concluding remarks}

We have also investigated groups with socle $P\Omega_{22}^\epsilon(q)$ or $P\Omega_{24}^\epsilon(q)$ using the techniques described in Section \ref{section:almostsimple}. We found that such groups do not act line-transitively on a finite linear space and we are therefore able to record the following result:

\begin{corollary}\label{corollary:final}
Let $G$ be a group which acts transitively on the set of lines of a linear space $\spaceS$. Suppose that $G$ has socle $L$ a simple classical group of dimension $n>20$. Then $G$ acts transitively on the set of flags of $\spaceS$ and we have Example \ref{example:flagtransitive}.
\end{corollary}

No doubt the techniques of Section \ref{section:almostsimple} can be applied to further reduce the lower bound in Corollary \ref{corollary:final}.

\bibliographystyle{amsalpha}
\bibliography{paper}

\newcommand{\etalchar}[1]{$^{#1}$}
\providecommand{\bysame}{\leavevmode\hbox to3em{\hrulefill}\thinspace}
\providecommand{\MR}{\relax\ifhmode\unskip\space\fi MR }
\providecommand{\MRhref}[2]{%
  \href{http://www.ams.org/mathscinet-getitem?mr=#1}{#2}
}
\providecommand{\href}[2]{#2}
\begin{thebibliography}{BDD{\etalchar{+}}90}

\bibitem[BDD]{bdd2}
F.~Buekenhout, A.~Delandtsheer, and J.~Doyen, \emph{Finite linear spaces with
  flag-transitive sporadic groups}, Unpublished notes.

\bibitem[BDD88]{bdd}
\bysame, \emph{Finite linear spaces with flag-transitive groups}, J. Combin.
  Theory, Series A \textbf{49} (1988), 268--293.

\bibitem[BDD{\etalchar{+}}90]{ftclass}
F.~Buekenhout, A.~Delandtsheer, J.~Doyen, P.~Kleidman, M.~Liebeck, and J.~Saxl,
  \emph{Linear spaces with flag-transitive automorphism groups}, Geom. Ded.
  \textbf{36} (1990), 89--94.

\bibitem[Blo67]{block}
R.E. Block, \emph{On the orbits of collineation groups}, Math. Zeitschrift
  \textbf{96} (1967), 33--49.

\bibitem[CKS76]{cks}
Charles~W. Curtis, William~M. Kantor, and Gary~M. Seitz, \emph{The
  {$2$}-transitive permutation representations of the finite {C}hevalley
  groups}, Trans. Amer. Math. Soc. \textbf{218} (1976), 1--59.

\bibitem[CNP03]{cnp}
A.~Camina, P.~Neumann, and C.~Praeger, \emph{Alternating groups acting on
  linear spaces}, Proc. London Math. Soc.(3) \textbf{87} (2003), no.~1, 29--53.

\bibitem[CP93]{campraeg2}
Alan~R. Camina and Cheryl~E. Praeger, \emph{Line-transitive automorphism groups
  of linear spaces}, Bull. London Math. Soc. \textbf{25} (1993), 309--315.

\bibitem[CP01]{campraeger}
A.~Camina and C.~Praeger, \emph{Line-transitive, point quasiprimitive
  automorphism groups of finite linear spaces are affine or almost simple},
  Aequationes Math. \textbf{61} (2001), 221--232.

\bibitem[CS89]{camsiem}
A.~Camina and J.~Siemons, \emph{Block transitive automorphism groups of
  2-$(v,k,1)$ block designs}, J. Combin. Theory, Series A \textbf{51} (1989),
  268--276.

\bibitem[CS00]{camspiez}
A.~Camina and F.~Spiezia, \emph{Sporadic groups and automorphisms of linear
  spaces}, J. Combin. Designs \textbf{8} (2000), 353--362.

\bibitem[Dav87]{davies}
D.H. Davies, \emph{Automorphism of designs}, Ph.D. thesis, University of East
  Anglia, 1987.

\bibitem[DD89]{deldoy}
A.~Delandtsheer and J.~Doyen, \emph{Most block-transitive $t$-designs are
  point-primitive}, Geom. Ded. \textbf{29} (1989), 307--310.

\bibitem[Del86]{delandt3}
A.~Delandtsheer, \emph{Flag-transitive finite simple groups}, Arch. Math.
  \textbf{47} (1986), 395--400.

\bibitem[Del01]{delandt4}
\bysame, \emph{Finite flag-transitive linear spaces with alternating socle},
  Algebraic combinatorics and applications (G\"o\ss weinstein, 1999), Springer,
  Berlin, 2001, pp.~79--88.

\bibitem[Gila]{gill1}
Nick Gill, \emph{${PSL}(3,q)$ and line-transitive linear spaces}, Beitr\"age
  zur Algebra und Geometrie, To appear.

\bibitem[Gilb]{gill2}
\bysame, \emph{Transitive projective planes}, Advances in Geometry, To appear.

\bibitem[Gil06]{gill3}
\bysame, \emph{Linear spaces with significant characteristic prime},
  Innovations in Incidence Geometry \textbf{3} (2006), 109--119.

\bibitem[GLS94]{gorenstein}
Daniel Gorenstein, Richard Lyons, and Ronald Solomon, \emph{The classification
  of the finite simple groups, number 3}, Mathematical Surveys and Monographs,
  vol.~40, American Mathematical Society, 1994.

\bibitem[GMS03]{gms}
Robert~M. Guralnick, Peter M{\"{u}}ller, and Jan Saxl, \emph{The rational
  function analogue of a question of {S}chur and exceptionality of permutation
  representations}, Mem. Amer. Math. Soc. \textbf{162} (2003), no.~773, 1--79.

\bibitem[KL90]{kl}
P.~Kleidman and M.~Liebeck, \emph{The subgroup structure of the finite simple
  groups}, London Mathematical Society Lecture Note Series, vol. 129, Cambridge
  University Press, Cambridge, 1990.

\bibitem[Kle90]{kleidman4}
Peter~B. Kleidman, \emph{The finite flag-transitive linear spaces with an
  exceptional automorphism group}, Finite geometries and combinatorial designs
  (Lincoln, NE, 1987), Contemp. Math., vol. 111, Amer. Math. Soc., Providence,
  RI, 1990, pp.~117--136.

\bibitem[Lie85]{liebeck3}
Martin~W. Liebeck, \emph{On the orders of maximal subgroups of the finite
  classical groups}, Proc. London Math. Soc. (3) \textbf{50} (1985), no.~3,
  426--446.

\bibitem[Liu01]{liug2even}
Weijun Liu, \emph{The {C}hevalley groups {$G\sb 2(2\sp n)$} and
  {$2\text{-}(v,k,1)$} designs}, Algebra Colloq. \textbf{8} (2001), no.~4,
  471--480.

\bibitem[Liu03a]{liug2odd}
\bysame, \emph{The {C}hevalley groups {$G\sb 2(q)$} with {$q$} odd and
  {$2\text{-}(v,k,1)$} designs}, European J. Combin. \textbf{24} (2003), no.~3,
  331--346.

\bibitem[Liu03b]{liupsu3even}
\bysame, \emph{Finite linear spaces admitting a projective group {${\rm
  PSU}(3,q)$} with {$q$} even}, Linear Algebra Appl. \textbf{374} (2003),
  291--305.

\bibitem[Liu03c]{liupsl2}
\bysame, \emph{Finite linear spaces admitting a two-dimensional projective
  linear group}, J. Combin. Theory Ser. A \textbf{103} (2003), no.~2, 209--222.

\bibitem[Liu03d]{liu3d4}
\bysame, \emph{{S}teinberg triality groups acting on 2-{$(v,k,1)$} designs},
  Sci. China Ser. A \textbf{46} (2003), no.~6, 872--883.

\bibitem[LLG06]{liuree2}
Weijun Liu, Shangzhao Li, and Luozheng Gong, \emph{Almost simple groups with
  socle {${\rm Ree}(q)$} acting on finite linear spaces}, European J. Combin.
  \textbf{27} (2006), no.~6, 788--800.

\bibitem[LLM01]{liusuzuki}
Weijun Liu, Huiling Li, and Chungui Ma, \emph{{S}uzuki groups and
  {$2$}-{$(v,k,1)$} designs}, European J. Combin. \textbf{22} (2001), no.~4,
  513--519.

\bibitem[LZLF04]{liuree}
Weijun Liu, Shenglin Zhou, Huiling Li, and Xingui Fang, \emph{Finite linear
  spaces admitting a {R}ee simple group}, European J. Combin. \textbf{25}
  (2004), no.~3, 311--325.

\bibitem[Sax02]{saxl}
Jan Saxl, \emph{On finite linear spaces with almost simple flag-transitive
  automorphism groups}, J. Combin. Theory A \textbf{100} (2002), 322--348.

\bibitem[Spi97]{fed}
F.~Spiezia, \emph{Simple groups and automorphisms of linear spaces}, {PhD}
  thesis, University of East Anglia, 1997.

\bibitem[Zho02]{zhouree2}
Shenglin Zhou, \emph{Block primitive {$2$}-{$(v,k,1)$} designs admitting a
  {R}ee simple group}, European J. Combin. \textbf{23} (2002), no.~8,
  1085--1090.

\bibitem[Zho05]{Zhouree}
\bysame, \emph{Block primitive 2-{$(v,k,1)$} designs admitting a {R}ee group of
  characteristic two}, Des. Codes Cryptogr. \textbf{36} (2005), no.~2,
  159--169.

\bibitem[ZLL00]{Zhouliu}
Shenglin Zhou, Huiling Li, and Weijun Liu, \emph{The {R}ee groups {${}\sp 2G\sb
  2(q)$} and {$2-(v,k,1)$} block designs}, Discrete Math. \textbf{224} (2000),
  no.~1-3, 251--258.

\end{thebibliography}

\end{document}